\documentclass[11pt]{article}
\usepackage{amsfonts}
\usepackage{latexsym,amsmath}
\usepackage{amssymb,array}
\usepackage{calc}
\RequirePackage[dvips]{graphicx}
\usepackage{graphics}
\usepackage[colorlinks=true]{hyperref}
\hypersetup{urlcolor=blue, citecolor=red}
\makeatletter \oddsidemargin 0in \evensidemargin 0in \textwidth
16cm\textheight 20cm 
\begin{document}
\newcommand{\la}{\lambda}
\newcommand{\eq}{\Leftrightarrow}
\newcommand{\mf}{\mathbf}
\newcommand{\ri}{\Rightarrow}
\newtheorem{t1}{Theorem}[section]
\newtheorem{d1}{Definition}[section]
\newtheorem{n1}{Notation}[section]
\newtheorem{c1}{Corollary}[section]
\newtheorem{l1}{Lemma}[section]
\newtheorem{r1}{Remark}[section]
\newtheorem{e1}{Counterexample}[section]
\newtheorem{re1}{Result}[section]
\newtheorem{p1}{Proposition}[section]
\newtheorem{cn1}{Conclusion}[section]
\renewcommand{\theequation}{\thesection.\arabic{equation}}
\pagenumbering{arabic}
\title {Ordering properties of sample minimum from Kumaraswamy-G random variables}
\author{Amarjit Kundu\\Department of
Mathematics\\
Santipur College\\ West Bengal, India
\and
Shovan Chowdhury\footnote{Corresponding
author e-mail: shovanc@iimk.ac.in; meetshovan@gmail.com}\;\\Quantitative Methods and Operations Management Area\\Indian Institute of Management, Kozhikode\\Kerala, India.}
\maketitle
\begin{abstract}
In this paper we compare the minimums of two independent and heterogeneous samples each following Kumaraswamy-G distribution with the same and the different parent distribution functions. The comparisons are carried out with respect to usual stochastic ordering and hazard rate ordering with majorized shape parameters of the distributions. The likelihood ratio ordering between the minimum order statistics is established for heterogeneous multiple outlier Kumaraswamy-G random variables with the same parent distribution function. 
\end{abstract}
\section{Introduction}
\setcounter{equation}{0}
\hspace*{0.3in} 
The paper by Kumaraswamy~\cite{kum} proposed a new two-parameter probability distribution on (0,1) with hydrological applications. The Kumaraswamy's distribution ($Kw$ distribution) does not seem to be popular in the statistical literature and has seen only limited use and development in the hydrological and related literatures (see Sundar and Subbiah~\cite{su}, Fletcher and Ponnambalam~\cite{fl}, Seifi \emph{et al.}~\cite{se} and Ganji \emph{et al.}~\cite{ga}). A recent paper by Jones~\cite{jo} explored the background and genesis of the $Kw$ distribution and discussed it's similarities to the beta distribution along with a number of advantages in terms of tractability. 
A random variable $X$ is said to have $Kw$ distribution with parameters ($\alpha,\beta$), written as $Kw$($\alpha,\beta$), if the cumulative distribution function (cdf) of $X$ is given by
\begin{equation*}\label{e0}
K(x)=1-\left(1-x^\alpha\right)^\beta,\;\;0<x<1,\;\alpha>0,\;\beta>0,
\end{equation*}
where $\alpha$ and $\beta$ are the shape parameters. Generalizing this distribution, Cordeiro and de Castro~\cite{co} have proposed a new family of generalized distributions, called Kumaraswamy generalized family of distributions (called $Kw$-G distribution). For a random variable $X$ with cdf $F(x)$, the distribution function $G(x)$ of the $Kw$-G random variable is defined by
\begin{equation}\label{e1}
G(x)=1-\left(1-F^\alpha(x)\right)^\beta,\;\;x\in \Re,\;\alpha>0,\;\beta>0.
\end{equation}
The $Kw$-G distribution, written as $Kw$-G($\alpha,\beta,F$), is shown to be used for the censored data quite effectively. Moreover, this distribution has the ability to fit skewed data better than any existing distributions. Each of the $Kw$-G distributions can be obtained from a specified parent cdf $F$, e.g. the $Kw$-normal ($Kw$-N) distribution is obtained by taking $F(x)$ as the cdf of the normal distribution. The $Kw$-Weibull ($Kw$-W), $Kw$-gamma ($Kw$-Ga) and $Kw$-Gumbel ($Kw$-Gu) distributions can be obtained similarly by taking $F(x)$ as the cdf of the Weibull, gamma and Gumbel distributions, respectively, among several others. Various properties of these distributions are discussed in the literature without attention to the stochastic properties of their order statistics.   
\\\hspace*{0.3in} Order statistics have a prominent role in reliability theory, life testing, actuarial science, auction theory, hydrology and many other related and unrelated areas. If $X_{1:n}\leq X_{2:n}\leq\ldots\leq X_{n:n}$ denote the order statistics corresponding to the random variables $X_1, X_2,\ldots,X_n$, then the sample minimum and sample maximum correspond to the smallest and the largest order statistics $X_{1:n}$ and $X_{n:n}$ respectively. For properties of order statistics for independent and non-identically distributed random variables, one may refer to David and Nagaraja~\cite{dn11}. The results of stochastic comparisons of the order statistics (largely on the smallest and the largest order statistics) can be seen in Dykstra \emph{et al.}~\cite{dkr11}, Fang and Zhang~(\cite{fz1},\cite{fz2}), Zhao and Balakrishnan~\cite{zb11.2}, Fang and Balakrishnan~\cite{fb}, Li and Li~\cite{li}, Torrado and Kochar~\cite{tr11}, Kundu \emph{et al.}~\cite{kun1}, Kundu and Chowdhury~\cite{kun2}, Chowdhury and Kundu~\cite{kun3} and the references there in for a variety of parametric models.
\\\hspace*{0.3in} In this paper our main aim is to compare minimums of two independent heterogenous samples from $Kw$-G random variables with both common ($F$) and different ($F_1$ and $F_2$) homogenous parent cdf. The comparison is carried out in terms of hazard rate order and likelihood ratio  order through majorization of the shape parameters. 
\\\hspace*{0.3 in} The organization of the paper is as follows. In Section 2, we have given the required definitions and some useful lemmas which are used throughout the paper. The results related to hazard rate ordering between two smallest order statistics from  two different $Kw$-G distributions having same parent distribution $F$, under the majorization of the shape parameters, are given in section~3. Here we have also shown that there exists likelihood ratio ordering between the smallest order statistics under certain arrangements of the parameters for multiple-outlier $Kw$-G model with the common homogenous parent cdf $F.$ Section~4 deals with various ordering related results between two smallest order statistics from two $Kw$-G distributions with different parent distributions. Finally, Section~5 concludes the paper. 
\\\hspace*{0.3 in}Throughout the paper, the word increasing (resp. decreasing) and nondecreasing (resp. nonincreasing) are used interchangeably, and $\Re$ denotes the set of real numbers $\{x:-\infty<x<\infty\}$. We also 
write $a\stackrel{sign}{=}b$ to mean that $a$ and $b$ have the same sign. For any differentiable function $k(\cdot)$,
we write $k'(t)$ to denote the first derivative of $k(t)$ with respect to $t$. 
\section{Preliminaries}
\hspace*{0.3 in} For two absolutely continuous random variables $X$ and $Y$ with distribution functions $F\left(\cdot\right)$ and $G\left(\cdot\right)$, survival functions $\overline F\left(\cdot\right)$ and $\overline G\left(\cdot\right)$, density functions $f\left(\cdot\right)$ and $g\left(\cdot\right)$ and hazard rate functions $r\left(\cdot\right)$ and $s\left(\cdot\right)$ respectively, $X$ is said to be smaller than $Y$ in $i)$ {\it likelihood ratio order} (denoted as $X\leq_{lr}Y$), if, for all $t$, $\frac{g(t)}{f(t)}$ increases in $t$, $ii)$ {\it hazard rate order} (denoted as $X\leq_{hr}Y$), if, for all $t$, $\frac{\overline G(t)}{\overline F(t)}$ increases in $t$ or equivalently $r(t)\geq s(t)$, and $iii)$ {\it usual stochastic order} (denoted as $X\leq_{st}Y$), if $F(t)\ge G(t)$ for all $t$. For more on different stochastic orders, see Shaked and Shanthikumar \cite{shak1}.
\\\hspace*{0.3 in} The notion of majorization (Marshall et al. [5]) is essential for the understanding of the
stochastic inequalities for comparing order statistics. Let $I^n$ be an $n$-dimensional Euclidean space where $I\subseteq\Re$. Further, for any two real vectors $\mathbf{x}=(x_1,x_2,\dots,x_n)\in I^n$ and $\mathbf{y}=(y_1,y_2,\dots,y_n)\in I^n$, write $x_{(1)}\le x_{(2)}\le\cdots\le x_{(n)}$ and $y_{(1)}\le y_{(2)}\le\cdots\le y_{(n)}$ as the increasing arrangements of the components of the vectors $\mathbf{x}$ and $\mathbf{y}$ respectively. The following definitions may be found in Marshall \emph{et al.} \cite{Maol}.
\begin{d1}
\begin{enumerate}
\item[i)] The vector $\mathbf{x} $ is said to majorize the vector $\mathbf{y} $ (written as $\mathbf{x}\stackrel{m}{\succeq}\mathbf{y}$) if
$$\sum_{i=1}^j x_{(i)}\le\sum_{i=1}^j y_{(i)},\;j=1,\;2,\;\ldots, n-1,\;\;and \;\;\sum_{i=1}^n x_{(i)}=\sum_{i=1}^n y_{(i)}.$$
\item [ii)] The vector $\mathbf{x}$ is said to weakly supermajorize the vector $\mathbf{y}$
 (written as $\mathbf{x}\stackrel{\rm w}{\succeq} \mathbf{y}$) if
  $$\sum\limits_{i=1}^j x_{(i)}\leq \sum\limits_{i=1}^j y_{(i)}\quad \text{for}\;j=1,2,\dots,n.$$
 \item [iii)] The vector $\mathbf{x}$ is said to weakly submajorize the vector $\mathbf{y}$
 (written as $\mathbf{x}\;{\succeq}_{\rm w} \;\mathbf{y}$) if
  $$\sum\limits_{i=j}^n x_{(i)}\geq \sum\limits_{i=j}^n y_{(i)}\quad \text{for}\;j=1,2,\dots,n.$$
 \end{enumerate}
\end{d1}
 Next we present some useful lemmas which will be used in the next section to prove our main results. The proof of the Lemmas \ref{l4}-\ref{l7} are straight forward and may be provided on request. The following lemma can be found in Marshall \emph{et al.} (\cite{Maol}, p. 87) where the parenthetical statements are not given.   
\begin{l1}\label{l2}
Let $\varphi: I^n\rightarrow \Re$. Then $$(a_1,a_2,\dots,a_n)\succeq_{w} (b_1,b_2,\dots,b_n)\; \text{implies} \;\varphi(a_1,a_2,\dots,a_n)\geq (\text{resp. }\leq)\; \varphi(b_1,b_2,\dots,b_n)$$ if, and only if, $\varphi$ is increasing (resp. decreasing) and Schur-convex (resp. Schur-concave) on $I^n$. Similarly,
 $$(a_1,a_2,\dots,a_n)\stackrel{w}\succeq (b_1,b_2,\dots,b_n)\; \text{implies} \;\varphi(a_1,a_2,\dots,a_n)\geq (\text{resp. }\leq)\;\varphi(b_1,b_2,\dots,b_n)$$ 
if, and only if, $\varphi$ is decreasing (resp. increasing) and Schur-convex (resp. Schur-concave) on $I^n$.
 \end {l1}
\begin{l1}\label{l3}
For $s,t>0$, and for any cdf $F(x),$ the function $\phi(s,t,x)=\frac{stF^{s}(x)}{1-F^{s}(x)}$ is decreasing in $s.$
\end{l1}
{\bf Proof:} Differentiating $\phi(s,t, x)$ partially with respect to $s$, we get $$\frac{\partial}{\partial s}\phi(s,t,x)=\frac{tF^{s}(x)\left(1-F^{s}(x)+\log F^{s}(x)\right)}{\left(1-F^{s}(x)\right)^2}.$$ Now, for all $x>0$, as $\log x \le x-1$, which implies that $\log F^{s}(x)\le F^{s}(x)-1$, then $\phi(s,t,x)$ is decreasing in $s$.
\begin{l1}\label{l4}
For $s>0$, and for any cdf $F(x),$ the function $\phi_1(s,x)=1+\frac{s\log F(x)}{1-F^{s}(x)}$ is decreasing in $s$, and consequently $\phi_2(s,x)=\frac{s}{1-F^{s}(x)}$ is increasing $s$.
\end{l1}
\begin{l1}\label{l6}
For $s>0$, and for any cdf $F(x),$ $\frac{\partial}{\partial s}\phi(s,t,x)=\frac{\phi(s,t, x)\phi_1(s,x)}{s}$ is decreasing in $s$.
\end{l1}
\begin{l1}\label{l7}
For $s>0$, and for any cdf $F(x),$ the function $\phi_3(s,x)=\frac{sF^{s}(x)\left(1-F^{s}(x)+sF^{s}(x)\log F(x)\right)}{\left(1-F^{s}(x)\right)^3}$ is decreasing in $s.$
\end{l1}
\begin{r1}\label{r2}
As  by Lemma \ref{l3}, $\phi(s,t, x)$ is decreasing in $s$ and $\phi(s,t,x) \ge 0$, then for all $s,x\ge 0$, $ \phi_1(s,x)\le 0$.
\end{r1}
\begin{n1}
Let us introduce the following notations which will be used in all the upcoming theorems.
\begin{enumerate}
\item[i)] $\mathcal{D}_{+}=\left\{\left(x_{1},x_2,\ldots,x_{n}\right):x_{1}\geq x_2\geq\ldots\geq x_{n}> 0\right\}$.
\item[ii)] $\mathcal{E}_{+}=\left\{\left(x_{1},x_2,\ldots,x_{n}\right):0< x_{1}\leq x_2\leq\ldots\leq x_{n}\right\}$.
\end{enumerate}
\end{n1}
\section{Results when Kw-G's have same parent distribution}
\setcounter{equation}{0}
Let $X$ be a random variable with continuous distribution function $F(\cdot)$ and density function $f(\cdot)$. Suppose that $U_i\sim$ $Kw$-G$\left(\alpha_i,\beta_i,F\right)$ and $V_i\sim$ $Kw$-G$\left(\gamma_i,\delta_i,F\right)$ ($i=1,2,\ldots,n$) be two sets of $n$ independent random variables where the parent cdf $F$ is homogenous and common to both the sets of random variables. Also suppose that $\overline{G}_{1:n}\left(\cdot\right)$ and $\overline{H}_{1:n}\left(\cdot\right)$ be the survival functions of $U_{1:n}$ and $V_{1:n}$ respectively. Then, for all $x\ge 0$, 
\begin{equation*}
\overline{G}_{1:n}\left(x\right)=\prod_{i=1}^n \left(1-F^{\alpha_i}(x)\right)^{\beta_i},
\end{equation*}
and
\begin{equation*}
\overline{H}_{1:n}\left(x\right)=\prod_{i=1}^n \left(1-F^{\gamma_i}(x)\right)^{\delta_i}.
\end{equation*}
Again, if $r_{1:n}(\cdot)$ and $s_{1:n}(\cdot)$ are the hazard rate functions of $U_{1:n}$ and $V_{1:n}$ respectively, then
\begin{equation}\label{e11}
r_{1:n}\left(x\right)=\sum_{i=1}^n\frac{\alpha_i\beta_i F^{\alpha_i-1}(x)f(x)}{1-F^{\alpha_i}(x)},
\end{equation} 
and
\begin{equation}\label{e21}
s_{1:n}\left(x\right)=\sum_{i=1}^n\frac{\gamma_i\delta_i F^{\gamma_i-1}(x)f(x)}{1-F^{\gamma_i}(x)}.
\end{equation} 
Let $\mbox{\boldmath $\alpha$}=\left(\alpha_1, \alpha_2, \ldots, \alpha_n\right),\ \mbox{\boldmath $\beta$}=\left(\beta_1, \beta_2, \ldots, \beta_n\right),\ \mbox{\boldmath $\gamma$}=\left(\gamma_1, \gamma_2, \ldots, \gamma_n\right)$ and $\mbox{\boldmath $\delta$}=\left(\delta_1, \delta_2, \ldots, \delta_n\right)\in I^n$. The following two theorems show that under certain conditions on parameters, there exists hazard rate ordering between $U_{1:n}$ and $V_{1:n}$. 
\begin{t1}\label{th7}
For $i=1,2,\ldots, n$, let $U_i$ and $V_i$ be two sets of mutually independent random variables with $U_i\sim$ $Kw$-G$\left(\alpha_i,\beta_i,F\right)$ and $V_i\sim$ $Kw$-G$\left(\gamma_i,\beta_i,F\right)$. Further, suppose that $\mbox{\boldmath $\alpha$},\mbox{\boldmath $\gamma$}\in \mathcal{D}_+$ and $\mbox{\boldmath $\beta$}\in \mathcal{E}_+$. Then, $\mbox{\boldmath $\alpha$}\stackrel{w}{\succeq}\mbox{\boldmath $\gamma$}\;\text{implies}\; X_{1:n}\le_{hr}Y_{1:n}.$
\end{t1}
{\bf Proof:} 
Let $a(\alpha_i)=\frac{\alpha_i F^{\alpha_i-1}(x)f(x)}{1-F^{\alpha_i}(x)}.$ Differentiating $a(\alpha_i)$ with respect to $\alpha_i$, we get $$a'(\alpha_i)=\frac{F^{\alpha_i-1}(x)f(x)\left(1-F^{\alpha_i}(x)+\log F^{\alpha_i}(x)\right)}{\left(1-F^{\alpha_i}(x)\right)^2}.$$  Following the proof of Lemma~\ref{l3} it can be proved that $a'(\alpha_i)\le 0.$ Hence, $a(\alpha_i)$ is decreasing in $\alpha_i.$ Again, differentiating $a'(\alpha_i)$ with respect to $\alpha_i$, we get  
\begin{equation}\label{e30}
a''(\alpha_i)=\frac{F^{\alpha_i-1}(x)f(x)\log F(x)}{\left(1-F^{\alpha_i}(x)\right)^3}b(\alpha_i),
\end{equation}
where, $b(\alpha_i)=2-2F^{\alpha_i}(x)+\log F^{\alpha_i}(x)+\alpha_iF^{\alpha_i}(x)\log F(x)$. Again, differentiating $b(\alpha_i)$ we get, 
\begin{equation}\label{e}
b'(\alpha_i)=\log F(x)c(\alpha_i),
\end{equation}
where $c(\alpha_i)=1-F^{\alpha_i(x)}+\alpha_i F^{\alpha_i}(x)\log F(x)$, which, on differentiation again with respect to $\alpha_i$ gives 
\begin{equation*}
c'(\alpha_i)=\alpha_i F^{\alpha_i}(x)\left(\log F(x)\right)^2>0.
\end{equation*}
Thus $c(\alpha_i)$ is increasing in $\alpha_i$ with $c(\alpha_i)=0$ at $\alpha_i=0.$ Therefore, for all $\alpha_i>0, c(\alpha_i)\ge 0$. So, equation (\ref{e}) gives $b(\alpha_i)$ is decreasing in $\alpha_i$ with $b(\alpha_i)=0$ at $\alpha_i=0.$ Hence, for all $\alpha_i>0,~b(\alpha_i) < 0$. So by equation (\ref{e30}), $a''(\alpha_i)\ge 0,$ giving that $a(\alpha_i)$ is convex in $\alpha_i.$ Thus by Theorem 3.1 b) (ii) of Kundu et al.~\cite{kun1} and Lemma \ref{l2} the result is proved.\hfill$\Box$   
\\\hspace*{0.2in} Theorem~\ref{th7} shows hr ordering between $U_{1:n}$ and $V_{1:n}$ when $\mbox{\boldmath $\alpha$}$ majorizes $\mbox{\boldmath $\gamma$}$ keeping the other parameters same. Now the question arises-what will happen if $\mbox{\boldmath $\beta$}$ majorizes $\mbox{\boldmath $\delta$}$ while the parameters $\mbox{\boldmath $\alpha$,\and\ $\gamma$}$ are equal? The theorem given below answers this question.
\begin{t1}\label{th8}
For $i=1,2,\ldots, n$, let $U_i$ and $V_i$ be two sets of mutually independent random variables with $U_i\sim$ $Kw$-G$\left(\alpha_i,\beta_i,F\right)$ and $V_i\sim $$Kw$-G$\left(\alpha_i,\delta_i,F\right)$. If $\mbox{\boldmath $\beta$}\stackrel{m}{\succeq}\mbox{\boldmath $\delta$}$ and 
\begin{enumerate}
\item[i)] $\mbox{\boldmath $\beta$},\mbox{\boldmath $\delta$}\in \mathcal{E}_+$ and $\mbox{\boldmath $\alpha$}\in \mathcal{D}_+$,  
then $U_{1:n}\le_{hr}V_{1:n}$;
\item[ii)] $\mbox{\boldmath $\beta$},\mbox{\boldmath $\delta$}\in \mathcal{D}_+$ and $\mbox{\boldmath $\alpha$}\in \mathcal{D}_+$,  
then $U_{1:n}\ge_{hr}V_{1:n}.$
\end{enumerate}
\end{t1}
{\bf Proof:}   
Equation~(\ref{e11}) can be written as
\begin{eqnarray*}
r(\mbox{\boldmath $\beta$},x)&=&\frac{f(x)}{F(x)}\sum_{i=1}^n\frac{\alpha_i\beta_iF^{\alpha_i}(x)}{1-F^{\alpha_i}(x)}\\&=&\frac{f(x)}{F(x)}\sum_{i=1}^n w_i \xi(\beta_i) (say),
\end{eqnarray*}
where $\xi(\beta_i)=\beta_i$ and $w_i=\frac{\alpha_iF^{\alpha_i}(x)}{1-F^{\alpha_i}(x)}$. Now, as by Lemma~\ref{l2}, $w_i$ is a decreasing function of $\alpha_i$, $\mbox{\boldmath $\alpha$} \in \mathcal{D}_+$ implies that $\mbox{\boldmath $w$} \in \mathcal{E}_+$. So, if $\mbox{\boldmath $\alpha$} \in \mathcal{D}_+$ i.e. if $\mbox{\boldmath $w$}\in \mathcal{E}_+$ and $\mbox{\boldmath $\beta$}\in \mathcal{E}_+$ ($\mbox{\boldmath $\beta$}\in \mathcal{D}_+$) then by Theorem 3.2 b)(i) (Theorem 3.1 b) (i)) of Kundu et al.~\cite{kun1} it can be proved that $r(\mbox{\boldmath $\beta$},x)$ is Schur convex on $\mathcal{E}_+$ (Schur concave on $\mathcal{D}_+$). This proves the result.\hfill$\Box$\\
The following theorem follows from Lemma \ref{l2} and Theorem \ref{th8}.
\begin{t1}
For $i=1,2,\ldots, n$, let $U_i$ and $V_i$ be two sets of mutually independent random variables with $U_i\sim$ $Kw$-G$\left(\alpha_i,\beta_i,F\right)$ and $V_i\sim$ $Kw$-G$\left(\alpha_i,\delta_i,F\right)$.  
\begin{enumerate}
\item[i)] If $\mbox{\boldmath $\beta$}\succeq_w\mbox{\boldmath $\delta$}$, $\mbox{\boldmath $\beta$},\mbox{\boldmath $\delta$}\in \mathcal{E}_+$ and $\mbox{\boldmath $\alpha$}\in \mathcal{D}_+$, then $U_{1:n}\le_{hr}V_{1:n}$;
\item[ii)] If $\mbox{\boldmath $\beta$}\stackrel{w}{\succeq}\mbox{\boldmath $\delta$}$, $\mbox{\boldmath $\beta$},\mbox{\boldmath $\delta$}\in \mathcal{D}_+$ and $\mbox{\boldmath $\alpha$}\in \mathcal{D}_+$, then $U_{1:n}\ge_{hr}V_{1:n}.$
\end{enumerate}
\end{t1}
The immediate question that can be raised- can the results of Theorem \ref{th7} and Theorem \ref{th8} be further extended to likelihood ratio (lr) ordering between $U_{1:n}$ and $V_{1:n}$ from hazard rate ordering? The next two counterexamples show that neither the result of Theorem \ref{th7} nor the result of Theorem \ref{th8} can be extended up to lr ordering for $n\ge 3$.\\
\hspace*{0.2in}The following counterexample shows that there does not exist lr ordering between $U_{1:n}$ and $V_{1:n}$ for $n\geq 3$ even if there exists majorization ordering between $\mbox{\boldmath $\alpha$}$ and $\mbox{\boldmath $\gamma$}$.
\begin{e1}\label{ce1}
Let $U_i\sim $$Kw$-G$\left(\alpha_i,\beta_i, F\right)$ and $V_i\sim$ $Kw$-G$\left(\gamma_i,\beta_i, F\right),$ $i=1,2,3$. Now, for $\left(\alpha_1,\alpha_2, \alpha_3\right)=\left(6.2,4.1,2\right)\in \mathcal{D}_+$ and $\left(\gamma_1,\gamma_2, \gamma_3\right)=\left(5.2,5.1,2\right)\in \mathcal{D}_+$, it is clear that $\mbox{\boldmath $\alpha$}\stackrel{m}{\succeq}\mbox{\boldmath $\gamma$}$. Now, if $\left(\beta_1, \beta_2, \beta_3\right)=\left(1,2,3\right)\in \mathcal{E}_+$ is taken, then from Figure 3.1 it can be concluded that $\frac{g_{1:3}(x)}{h_{1:3}(x)}$ is not monotone where $g_{1:3}$ and $h_{1:3}$ are the pdf's of the random variables $U_{1:3}$ and $V_{1:3}$ respectively. 
\begin{figure}[t]\centering
\includegraphics[height=7 cm]{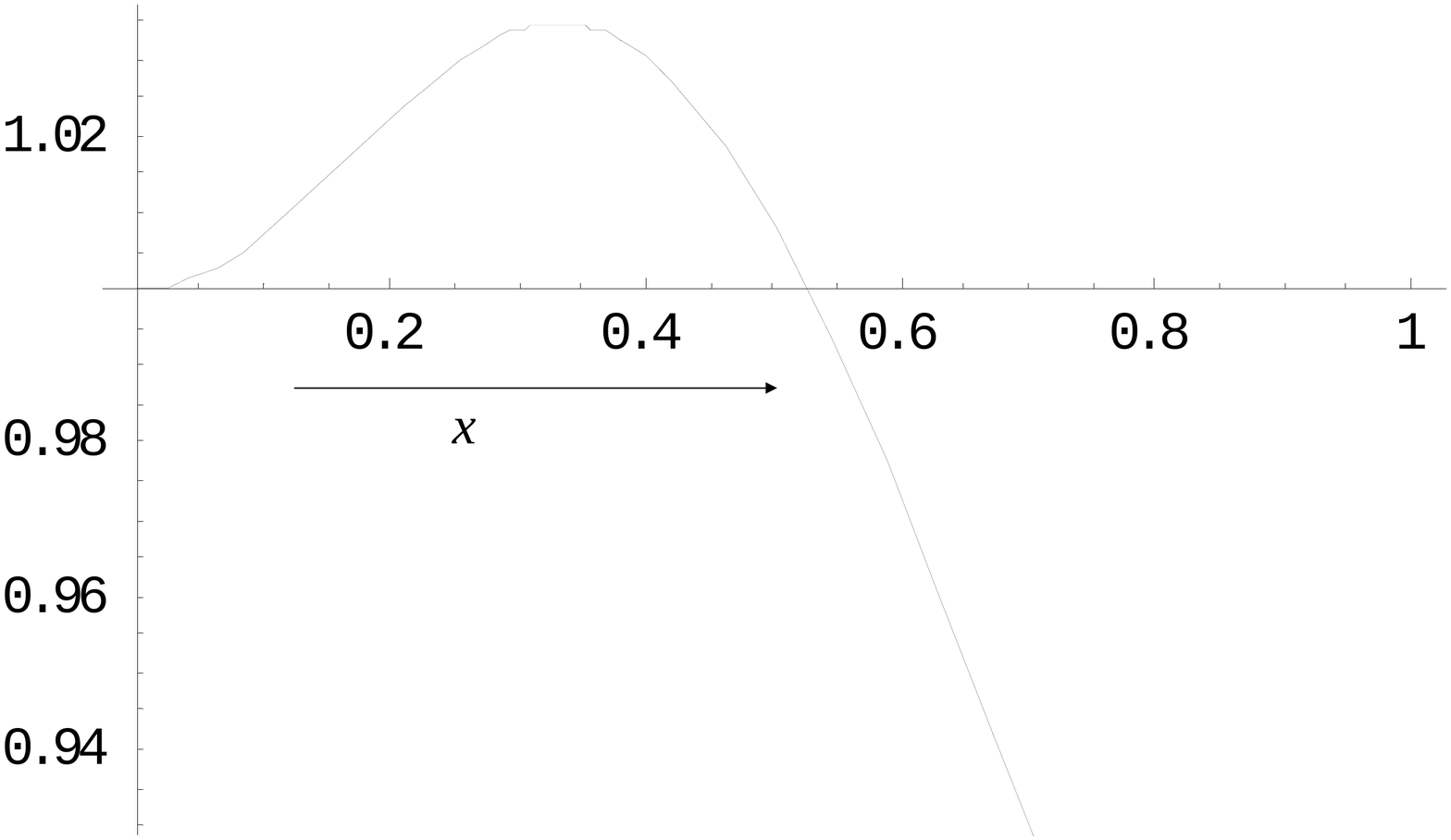}
\center{\it Figure 3.1: Graph of $\frac{g_{1:3}(x)}{h_{1:3}(x)}$} 
\end{figure}
\end{e1}
\hspace*{0.2in}The next counterexample shows that for $n\ge 3$ the result of Theorem \ref{th8} cannot be extended up to lr ordering.
\begin{e1}\label{ce2}
Let $U_i\sim $$Kw$-G$\left(\alpha_i,\beta_i, F\right)$ and $V_i\sim $$Kw$-G$\left(\alpha_i,\delta_i, F\right),$ $i=1,2,3$. Now, for $\left(\alpha_1,\alpha_2, \alpha_3\right)=\left(5,1,0.01\right)\in \mathcal{D}_+,$ $\left(\beta_1,\beta_2,\beta_3\right)=\left(0.005,0.004,0.001\right)\in \mathcal{D}_+$ and $\left(\delta_1,\delta_2,\delta_3\right)=\left(0.0045,0.0045,0.001\right)\in \mathcal{D}_+$, it is clear that $\mbox{\boldmath $\beta$}\stackrel{m}{\succeq}\mbox{\boldmath $\delta$}$. But, Figure 3.2 (i) shows that $\frac{g_{1:3}(x)}{h_{1:3}(x)}$ is not monotone. Again, for same $\mbox{\boldmath $\alpha$},$ if $\left(\beta_1,\beta_2,\beta_3\right)=\left(0.003,0.004,0.005\right)\in \mathcal{E}_+$ and $\left(\delta_1,\delta_2,\delta_3\right)=\left(0.0035,0.0035,0.005\right)\in \mathcal{E}_+$ are taken, then it can noticed that although $\mbox{\boldmath $\beta$}\stackrel{m}{\succeq}\mbox{\boldmath $\delta$},$ $\frac{g_{1:3}(x)}{h_{1:3}(x)}$ is non-monotone as evident from Figure 3.2 (ii). 
\begin{figure}[ht]
\centering
\begin{minipage}[b]{0.45\linewidth}
\includegraphics[height=7 cm]{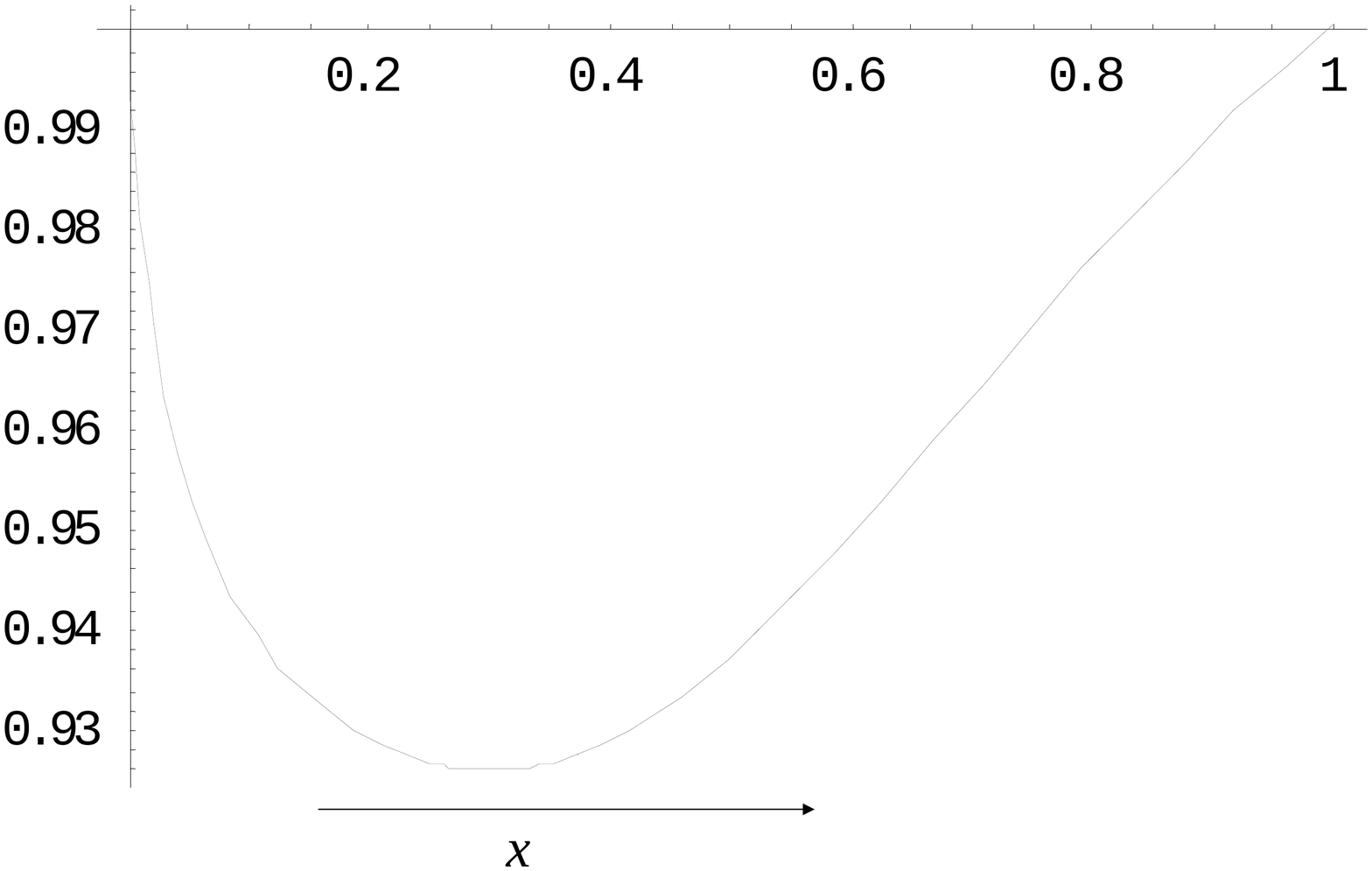}
\centering{$\left(i\right)$ Graph for $\mbox{\boldmath $\alpha$}, \mbox{\boldmath $\beta$}, \mbox{\boldmath $\delta$}\in \mathcal{D}_+$ }
\end{minipage}
\quad
\begin{minipage}[b]{0.45\linewidth}
\includegraphics[height=7 cm]{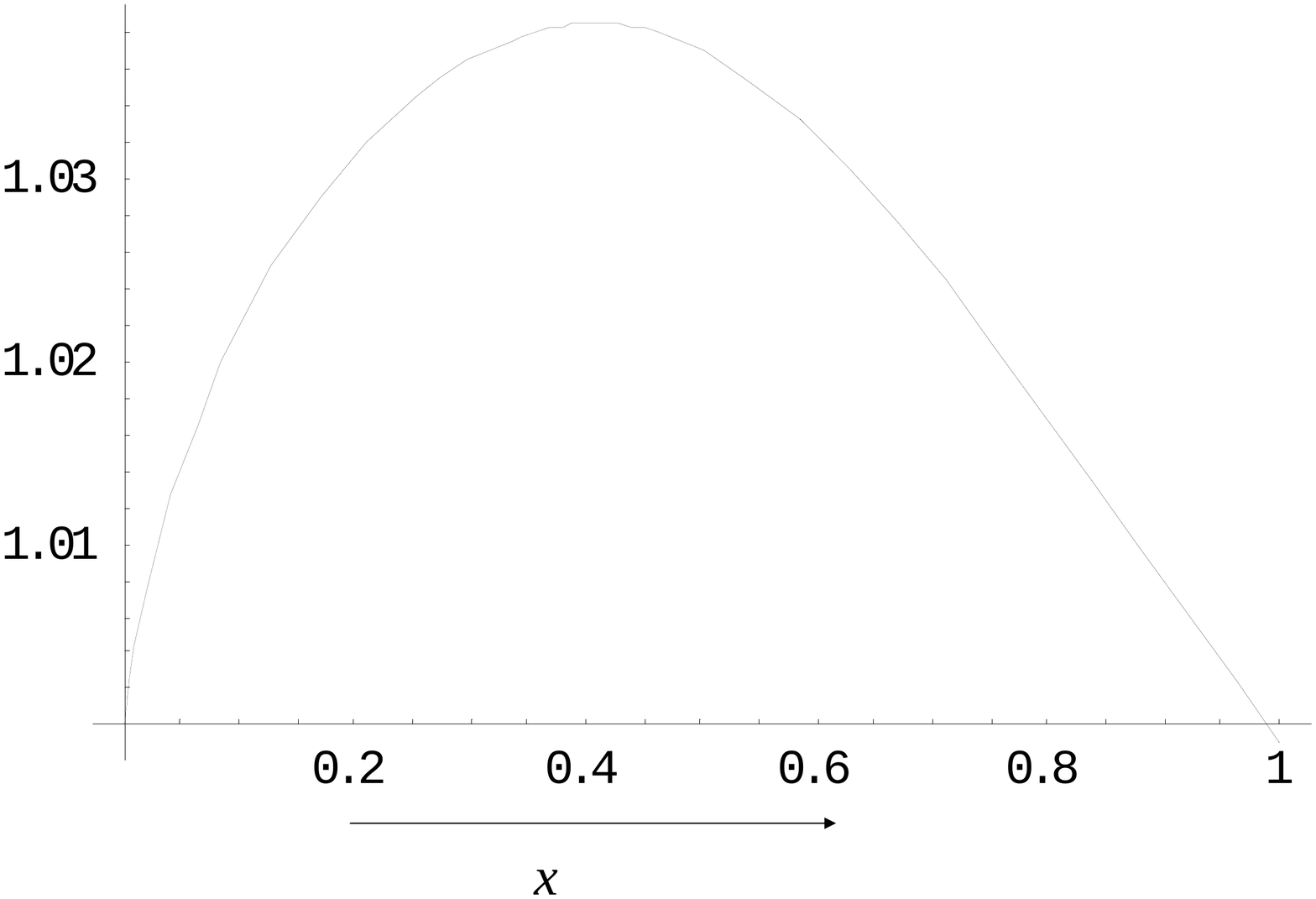}
\centering{$\left(ii\right)$ Graph for $\mbox{\boldmath $\alpha$}\in \mathcal{D}_+, \mbox{\boldmath $\beta$}, \mbox{\boldmath $\delta$}\in \mathcal{E}_+$}
\end{minipage}
\center{\it Figure 3.2: Graph of $\frac{g_{1:3}(x)}{h_{1:3}(x)}$} 
\end{figure}
\end{e1}
\hspace*{0.2in} Although there exists no lr ordering between $U_{1:n}$ and $V_{1:n}$ for $n\ge 3$, the following two theorems show that in case of multiple-outlier model lr ordering exists for any positive integer $n$.
\begin{t1}\label{th9}
 For $i=1,2,...,n$, let $U_i$ and $V_i$ be two sets of mutually independent random variables each following multiple-outlier $Kw$-G model such that $U_i\sim $$Kw$-G$\left(\alpha,\beta,F\right)$ and $V_i\sim$ $Kw$-G$\left(\gamma,\beta,F\right)$ for $i=1,2,\ldots,n_1$,
 $U_i\sim$ $Kw$-G$\left(\alpha^*,\beta^*, F\right)$ and $V_i\sim$ $Kw$-G$\left(\gamma^*,\beta^*, F\right)$ for $i=n_1+1,n_1+2,\ldots,n_1+n_2(=n)$. If $\alpha>\alpha^*, \gamma>\gamma^*$ and $\beta<\beta^*,$ and if $$(\underbrace{\alpha,\alpha,\ldots,\alpha,}_{n_1} \underbrace{\alpha^*,\alpha^*,\ldots,\alpha^*}_{n_2})\stackrel{m}{\succeq} (\underbrace{\gamma,\gamma,\ldots,\gamma,}_{n_1} \underbrace{\gamma^*,\gamma^*,\ldots,\gamma^*}_{n_2}),~then~U_{1:n}\le_{lr}V_{1:n}$$. 
\end{t1}
{\bf Proof:} 
In view of Theorem~\ref{th7} we need only to prove that $\frac{r_{1:n}(x)}{s_{1:n}(x)}$ is decreasing in $x.$ Now, 
\begin{equation*}
\frac{d}{dx}\left(\frac{r_{1:n}(x)}{s_{1:n}(x)}\right)\stackrel{sign}=\left[\sum_{k=1}^{n}\frac{\alpha_{k}^2\beta_kF^{\alpha_{k}}(x)}{(1-F^{\alpha_k}(x))^2}\right]\left[\sum_{k=1}^{n}\frac{\gamma_{k}\beta_kF^{\gamma_k}(x)}{1-F^{\gamma_k}(x)}\right]-\left[\sum_{k=1}^{n}\frac{\gamma_{k}^2\beta_kF^{\gamma_k}(x)}{(1-F^{\gamma_k}(x))^2}\right]\left[\sum_{k=1}^{n}\frac{\alpha_{k}\beta_kF^{\alpha_k}(x)}{1-F^{\alpha_k}(x)}\right], 
\end{equation*}
where $\alpha_k=\alpha, \beta_k=\beta, \gamma_k=\gamma$ for $k=1,2,\ldots, n_1$ and $\alpha_k=\alpha^*, \beta_k=\beta^*, \gamma_k=\gamma^*$ for $k=n_1+1,n_1+2,\ldots, n_1+n_2(=n).$
Thus, to show that $\frac{r_{1:n}(x)}{s_{1:n}(x)}$ is decreasing in $x,$ it is sufficient to show that
$$\Psi(\mbox{\boldmath $\alpha$},x)=\frac{\sum_{k=1}^{n}\frac{\alpha_{k}^2\beta_kF^{\alpha_k}(x)}{(1-F^{\alpha_k}(x))^2}}{\sum_{k=1}^{n}\frac{\alpha_k\beta_kF^{\alpha_k}(x)}{1-F^{\alpha_k}(x)}}=\frac{\sum_{k=1}^{n}\phi(\alpha_{k},\beta_k, x)\phi_2(\alpha_k,x)}{\sum_{k=1}^{n}\phi(\alpha_{k}, \beta_k, x)}$$
is Schur-concave in $\mbox{\boldmath $\alpha$},$ where $\phi(\alpha_k, \beta_k,x)$ and $\phi_2(\alpha_k,x)$ are same as defined in Lemma~\ref{l3} and Lemma~\ref{l4} respectively. Differentiating $\Psi$ partially with respect to $\alpha$ and $\alpha^*$ and using Lemma \ref{l6} we get,
\begin{equation*}
\begin{split}
\frac{\partial \Psi}{\partial \alpha}&\stackrel{sign}{=}\phi(\alpha,\beta,x)\left[n_1\phi(\alpha,\beta,x)\frac{\partial\phi_2(\alpha,x)}{\partial\alpha}+n_2\phi(\alpha^{*},\beta^{*}, x)\frac{\partial\phi_2(\alpha,x)}{\partial\alpha}\right.\\&\quad+\left.\frac{n_2\phi(\alpha^{*},\beta^{*},x)\phi_2(\alpha,x)}{\alpha}\left(\phi_1(\alpha,x)-\phi_1(\alpha^{*},x)\right)\right],
\end{split}
\end{equation*}
with a similar expression for $\frac{\partial \Psi}{\partial \alpha^{*}}$, where $\phi_1(\alpha,x)$ is same as defined in Lemma~\ref{l4}. \\
\hspace*{0.3 in} Now, three cases may arise:\\
$Case (i)$ $1\leq i<j\leq n_1.$ Here $\alpha_i=\alpha_j=\alpha$ and $\beta_i=\beta_j=\beta$, so that
$$\frac{\partial \Psi}{\partial \alpha_i}-\frac{\partial \Psi}{\partial \alpha_j}=\frac{\partial \Psi}{\partial \alpha}-\frac{\partial \Psi}{\partial \alpha}=0.$$
$Case (ii)$ If $n_1+1\leq i<j\leq n$, $i.e.$ if $\alpha_i=\alpha_j=\alpha^*$ and $\beta_i=\beta_j=\beta^*$, then
$$\frac{\partial \Psi}{\partial \alpha_i}-\frac{\partial \Psi}{\partial \alpha_j}=\frac{\partial \Psi}{\partial \alpha^*}-\frac{\partial \Psi}{\partial \alpha^*}=0.$$
$Case (iii)$ If $1\leq i\leq n_1$ and $n_1+1\leq j\leq n$, then $\alpha_i=\alpha$, $\beta_i=\beta$ and $\alpha_j=\alpha^*$, $\beta_j=\beta^*$. Then, 
\begin{eqnarray}\label{e7}
\begin{split}
\frac{\partial \Psi}{\partial \alpha_i}-\frac{\partial \Psi}{\partial \alpha_j}=&\frac{\partial \Psi}{\partial \alpha}-\frac{\partial \Psi}{\partial \alpha^*}\\
\stackrel{sign}{=}&\left[(n_1\phi(\alpha,\beta, x)+n_2\phi(\alpha^{*},\beta^*,x))\left(\phi(\alpha,\beta,x)\frac{\partial\phi_2(\alpha,x)}{\partial\alpha}-\phi(\alpha^{*},\beta^*,x)\frac{\partial\phi_2(\alpha^{*},x)}{\partial\alpha^{*}}\right)\right]\\&\quad+\left[\phi(\alpha,\beta,x)\phi(\alpha^{*},\beta^*,x)\left(\phi_2(\alpha,x)-\phi_2(\alpha^{*},x)\right)\left(\frac{n_1\phi_1(\alpha^*,x)}{\alpha^*}+\frac{n_2\phi_1(\alpha,x)}{\alpha}\right)\right].
\end{split}
\end{eqnarray}
Now, for all $\alpha,x\ge 0$, as $\phi_1(\alpha,x)<0$ (by Remark~\ref{r2}) and $\phi_2(\alpha,x)$ is increasing in $\alpha$ (by Lemma~\ref{l4}), and as $\alpha>\alpha^{*}$, then the second bracketed term of (\ref{e7}) is negative. Again, from Lemma~\ref{l7} it can be written that
\begin{eqnarray*}
\phi(\alpha,\beta,x)\frac{\partial\phi_2(\alpha,x)}{\partial\alpha}&=&\frac{\alpha\beta F^{\alpha}(x)\left(1-F^{\alpha}(x)+\alpha F^{\alpha}(x)\log F(x)\right)}{\left(1-F^{\alpha}(x)\right)^3}\\&<&\frac{\alpha^*\beta F^{\alpha^*}(x)\left(1-F^{\alpha^*}(x)+\alpha^* F^{\alpha^*}(x)\log F(x)\right)}{\left(1-F^{\alpha^*}(x)\right)^3}\\&<&\frac{\alpha^*\beta^* F^{\alpha^*}(x)\left(1-F^{\alpha^*}(x)+\alpha^* F^{\alpha^*}(x)\log F(x)\right)}{\left(1-F^{\alpha^*}(x)\right)^3}\\&=&\phi(\alpha^*,\beta^*,x)\frac{\partial\phi_2(\alpha^{*},x)}{\partial\alpha^{*}},
\end{eqnarray*}
where the second inequality follows from the facts that $\beta<\beta^{*}$, and for all $\alpha^*\ge 0$, $\frac{\alpha^*F^{\alpha^*}(x)}{1-F^{\alpha^*}(x)}>0,$ and $\frac{\partial\phi_2(\alpha,x)}{\partial\alpha}>0$ (by Lemma~\ref{l4}). So the first bracketed term of (\ref{e7}) is also negative. Thus, for all $i\le j$ it can be written that 
$$\frac{\partial \Psi}{\partial \alpha_i}-\frac{\partial \Psi}{\partial \alpha_j}\le 0, $$proving the result by Lemma 3.1 of  Kundu et al.~\cite{kun1}.\hfill$\Box$
\begin{t1}\label{th10}
 For $i=1,2,...,n$, let $U_i$ and $V_i$ be two sets of mutually independent random variables each following the multiple-outlier $Kw$-G model such that $U_i\sim $$Kw$-G$\left(\alpha,\beta,F\right)$ and $V_i\sim $$Kw$-G$\left(\alpha,\delta,F\right)$ for $i=1,2,\ldots,n_1$,
 $U_i\sim$ $Kw$-G$\left(\alpha^*,\beta^*,F\right)$ and $V_i\sim$ $Kw$-G$\left(\alpha^*,\delta^*,F\right)$ for $i=n_1+1,n_1+2,\ldots,n_1+n_2(=n)$. 
If $$(\underbrace{\beta,\beta,\ldots,\beta,}_{n_1} \underbrace{\beta^*,\beta^*,\ldots,\beta^*}_{n_2})\stackrel{m}{\succeq} (\underbrace{\delta,\delta,\ldots,\delta,}_{n_1} \underbrace{\delta^*,\delta^*,\ldots,\delta^*}_{n_2})$$ and
\begin{enumerate}
\item[i)] $\alpha>\alpha^*, \beta>\beta^*, \delta<\delta^*,$ then $U_{1:n}\ge_{lr}V_{1:n}$; 
\item[ii)] $\alpha>\alpha^*, \beta<\beta^*, \delta<\delta^*,$  then $U_{1:n}\le_{lr}V_{1:n}$. 
\end{enumerate}
\end{t1}
{\bf Proof:} 
In view of Theorem~\ref{th8} we need only to prove that $\frac{r_{1:n}(x)}{s_{1:n}(x)}$ is increasing in $x$ under conditions $i)$ and decreasing in $x$ under conditions $ii)$. Now, 
\begin{equation*}
\frac{d}{dx}\left(\frac{r_{1:n}(x)}{s_{1:n}(x)}\right)\stackrel{sign}=\left[\sum_{i=1}^{n}\frac{\alpha_{i}^2\beta_iF^{\alpha_i}(x)}{(1-F^{\alpha_i}(x))^2}\right]\left[\sum_{i=1}^{n}\frac{\alpha_{i}\delta_iF^{\alpha_i}(x)}{1-F^{\alpha_i}(x)}\right]-\left[\sum_{i=1}^{n}\frac{\alpha_{i}^2\delta_iF^{\alpha_i}(x)}{(1-F^{\alpha_i}(x))^2}\right]\left[\sum_{i=1}^{n}\frac{\alpha_{i}\beta_iF^{\alpha_i}(x)}{1-F^{\alpha_i}(x)}\right].
\end{equation*}
Thus, to show that $\frac{r_{1:n}(x)}{s_{1:n}(x)}$ is increasing (decreasing) in $x,$ it is sufficient to show that
$$\Psi(\mbox{\boldmath $\beta$},x)=\frac{\sum_{i=1}^{n}\frac{\alpha_{i}^2\beta_iF^{\alpha_i}(x)}{(1-F^{\alpha_i}(x))^2}}{\sum_{i=1}^{n}\frac{\alpha_i\beta_iF^{\alpha_i}(x)}{1-F^{\alpha_i}(x)}}$$
is Schur-convex (Schur-concave) in $\mbox{\boldmath $\beta$}.$ Now,
\begin{equation*}
\frac{\partial \Psi}{\partial \beta_i}\stackrel{sign}{=}\frac{\alpha_{i}^{2}F^{\alpha_i}(x)}{(1-F^{\alpha_i}(x))^2}\left[\frac{n_1\alpha\beta F^{\alpha}(x)}{1-F^{\alpha}(x)}+\frac{n_2\alpha^{*}\beta^{*}F^{\alpha^*}(x)}{1-F^{\alpha^*}(x)}\right]-\frac{\alpha_iF^{\alpha_i}(x)}{1-F^{\alpha_i}(x)}\left[\frac{n_1\alpha^2\beta F^{\alpha}(x)}{(1-F^{\alpha}(x))^2}+\frac{n_2\alpha^{*2}\beta^* F^{\alpha^*}(x)}{(1-F^{\alpha^*}(x))^2}\right].
\end{equation*}
\hspace*{0.3 in} Now, three cases may arise:\\
$Case\ i)$ If $1\leq i<j\leq n_1$, $i.e.$ if $\alpha_i=\alpha_j=\alpha$ and $\beta_i=\beta_j=\beta$, then $$\frac{\partial \Psi}{\partial \beta_i}-\frac{\partial \Psi}{\partial \beta_j}=\frac{\partial \Psi}{\partial \beta}-\frac{\partial \Psi}{\partial \beta}=0.$$ 
$Case\ ii)$ Again, if $n_1+1\leq i<j\leq n,$ then $\alpha_i=\alpha_j=\alpha^*$ and $\beta_i=\beta_j=\beta^*$ and correspondingly $$\frac{\partial \Psi}{\partial \beta_i}-\frac{\partial \Psi}{\partial \beta_j}=\frac{\partial \Psi}{\partial \beta^*}-\frac{\partial \Psi}{\partial \beta^*}=0.$$
$Case\ iii)$ Now, if $1\leq i\leq n_1$ and $n_1+1\leq j\leq n$, then $\alpha_i=\alpha$, $\beta_i=\beta$, and $\alpha_j=\alpha^*$ and $\beta_j=\beta^*$. So, 
\begin{eqnarray*}
\frac{\partial \Psi}{\partial \beta_i}-\frac{\partial \Psi}{\partial \beta_j}&=&\frac{\partial \Psi}{\partial \beta}-\frac{\partial \Psi}{\partial \beta^*}\\&=&\left(n_1\beta+n_2\beta^*\right)\frac{\alpha\alpha^*F^{\alpha+\alpha^*}(x)}{\left(1-F^{\alpha}(x)\right)\left(1-F^{\alpha^*}(x)\right)}\left[\frac{\alpha}{1-F^{\alpha}(x)}-\frac{\alpha^*}{1-F^{\alpha^*}(x)}\right].
\end{eqnarray*}
Now, if $\alpha>\alpha^{*}$, then from Lemma~\ref{l4} we have $\frac{\alpha}{1-F^{\alpha}(x)}-\frac{\alpha^*}{1-F^{\alpha^*}(x)}>0$.  Thus, for all $x\ge 0$ and for all $i\le j$ it can be written that
$$\frac{\partial \Psi}{\partial \beta_i}-\frac{\partial \Psi}{\partial \beta_j}\geq 0.$$  Hence, by Lemma 3.1 of Kundu et al.~\cite{kun1}, $\Psi$ is Schur-convex in $\mbox{\boldmath $\beta$}\in \mathcal{D}_+$. Again, by Lemma 3.2 of Kundu et al.~\cite{kun1}, $\Psi$ is Schur-concave in $\mbox{\boldmath $\beta$}\in \mathcal{E}_+$. This proves the result.\hfill$\Box$\\
\begin{r1}\label{r3}
Although counterxamples \ref{ce1} and \ref{ce2} showed that there exits no lr ordering between $U_{1:n}$ and $V_{1:n}$ for $n\ge 3$, the above two theorems show that the results are true for $n=2$.  
\end{r1}
\section{Results when Kw-G's have different parent distributions}
\setcounter{equation}{0}
In this section we generalize the previous model by taking two $Kw$-G random variables with different homogenous parent cdf's. Let $X_1$ and $X_2$ be two random variables with continuous distribution functions $F_1(\cdot)$ and $F_2(\cdot)$ and density functions $f_1(\cdot)$ and $f_2(\cdot)$ respectively. Also suppose that $U_i\sim$ $Kw$-G$\left(\alpha_i,\beta_i,F_1\right)$ and $V_i\sim$ $Kw$-G$\left(\gamma_i,\delta_i,F_2\right)$ ($i=1,2,\ldots,n$) be two sets of $n$ independent random variables. Therefore, for all $x\ge 0$
\begin{equation*}
\overline{G}_{1:n}\left(x\right)=\prod_{i=1}^n \left(1-F_1^{\alpha_i}(x)\right)^{\beta_i},
\end{equation*}
and
\begin{equation*}
\overline{H}_{1:n}\left(x\right)=\prod_{i=1}^n \left(1-F_2^{\gamma_i}(x)\right)^{\delta_i},
\end{equation*} 
represent survival functions of $U_{1:n}$ and $V_{1:n}$ respectively.\\
The next two theorems show that under certain conditions on parameters usual stochastic ordering between $X_1$ and $X_2$ implies the same between $U_{1:n}$ and $V_{1:n}$. 
\begin{t1}\label{th11}
For $i=1,2,\ldots, n$, let $U_i$ and $V_i$ be two sets of mutually independent random variables with $U_i\sim $$Kw$-G$\left(\alpha_i,\beta_i,F_1\right)$ and $V_i\sim$ $Kw$-G$\left(\gamma_i,\beta_i,F_2\right)$. Further, suppose that $\mbox{\boldmath $\alpha$},\mbox{\boldmath $\gamma$}\in \mathcal{D}_+$ and $\mbox{\boldmath $\beta$}\in \mathcal{E}_+$ and $\mbox{\boldmath $\alpha$}\stackrel{m}{\succeq}\mbox{\boldmath $\gamma$}$, then, $X_1\le_{st}X_2$ implies $U_{1:n}\le_{st}V_{1:n}.$
\end{t1}
{\bf Proof:} 
Let us consider another random variable $W_i$ such that $W_i\sim $$Kw$-G$\left(\gamma_i,\beta_i,F_1\right)$. If $\mbox{\boldmath $\alpha$},\mbox{\boldmath $\gamma$}\in \mathcal{D}_+$ and $\mbox{\boldmath $\beta$}\in \mathcal{E}_+,$ then by Theorem~\ref{th7}, it can be shown that   
$U_{1:n}\le_{hr}W_{1:n}$, which implies that $U_{1:n}\le_{st}W_{1:n}$. Thus, by definition of st ordering it can be written that   
\begin{equation}\label{4.1}
\prod_{i=1}^n \left(1-F_{1}^{\alpha_i}(x)\right)^{\beta_i}\le \prod_{i=1}^n \left(1-F_{1}^{\gamma_i}(x)\right)^{\beta_i}.
\end{equation}
Now, $X_1\le_{st}X_2$ implies $F_1(x)>F_2(x)$, which gives $1- F_{1}^{\gamma_i}(x)\le 1-F_{2}^{\gamma_i}(x)$ for all $\gamma_i\ge 0$. So from (\ref{4.1}), it can be written that $$\prod_{i=1}^n \left(1-F_{1}^{\alpha_i}(x)\right)^{\beta_i}\leq \prod_{i=1}^n \left(1-F_{1}^{\gamma_i}(x)\right)^{\beta_i}\leq \prod_{i=1}^n \left(1-F_{2}^{\gamma_i}(x)\right)^{\beta_i},$$ giving $U_{1:n}\le_{st}V_{1:n}.$   \hfill$\Box$  
\begin{t1}\label{th13}
For $i=1,2,\ldots, n$, let $U_i$ and $V_i$ be two sets of mutually independent random variables with $U_i\sim $$Kw$-G$\left(\alpha_i,\beta_i,F_1\right)$ and $V_i\sim$ $Kw$-G$\left(\alpha_i,\delta_i,F_2\right)$. If $\mbox{\boldmath $\beta$}\stackrel{m}{\succeq}\mbox{\boldmath $\delta$},$ $\mbox{\boldmath $\alpha$}\in \mathcal{D}_+$ and
\begin{enumerate}
\item[i)] $\mbox{\boldmath $\beta$},\mbox{\boldmath $\delta$}\in \mathcal{E}_+$, then $X_1\le_{st}X_2$ implies $U_{1:n}\le_{st}V_{1:n}$;
\item[ii)] $\mbox{\boldmath $\beta$},\mbox{\boldmath $\delta$}\in \mathcal{D}_+$, then $X_1\ge_{st}X_2$ implies $U_{1:n}\ge_{st}V_{1:n}.$
\end{enumerate}
\end{t1}
{\bf Proof:}
Considering $W_i\sim $$Kw$-G$\left(\alpha_i,\delta_i,F_1\right)$, and using the same logic as of Theorem~\ref{th11}, the theorem can be proved with the help of Theorem~\ref{th8}.\hfill$\Box$\\
Now the question arises, whether the results of Theorem~\ref{th11} and Theorem~\ref{th13} can be extended to hr ordering? The next two theorems answers this question.
\begin{t1}\label{th12}
For $i=1,2,\ldots, n$, let $U_i$ and $V_i$ be two sets of mutually independent random variables with $U_i\sim $$Kw$-G$\left(\alpha_i,\beta_i,F_1\right)$ and $V_i\sim$ $Kw$-G$\left(\gamma_i,\beta_i,F_2\right)$. Further, for any real number $s>0$ and for $i=1,2$, suppose that $X_i^s$ be random variables having distribution functions $F_{i}^s(x)$. Now, if $\mbox{\boldmath $\alpha$},\mbox{\boldmath $\gamma$}\in \mathcal{D}_+$, $\mbox{\boldmath $\beta$}\in \mathcal{E}_+$ and for all $s>0$ $X_1^{s}\le_{hr}X_2^{s}$, then, $\mbox{\boldmath $\alpha$}\stackrel{m}{\succeq}\mbox{\boldmath $\gamma$}\;\text{implies}\; U_{1:n}\le_{hr}V_{1:n}.$
\end{t1}
{\bf Proof:}  
Let us consider the random variable $W_i$ as defined in Theorem~\ref{th11}. As $\mbox{\boldmath $\alpha$},\mbox{\boldmath $\gamma$}\in \mathcal{D}_+$ and $\mbox{\boldmath $\beta$}\in \mathcal{E}_+,$ by Theorem~\ref{th7}, it can be written that 
\begin{equation}\label{e61}
\frac{\prod_{i=1}^n \left(1-F_{1}^{\gamma_{i}}(x)\right)^{\beta_i}}{\prod_{i=1}^n \left(1-F_{1}^{\alpha_{i}}(x)\right)^{\beta_i}}~is~increasing~in~x.
\end{equation}  
Again, for all $s\ge 0$, $X_1^{s}\le_{hr}X_2^{s}$ implies that $\frac{1-F_{2}^s(x)}{1-F_{1}^s(x)}$ is increasing in $x$. Thus, for all $\gamma_i>0,$ $\frac{1-F_{2}^{\gamma_i}(x)}{1-F_{1}^{\gamma_i}(x)}$ is increasing in $x$ implies that $\sum_{i=1}^n \beta_i \log\left(1-F_{2}^{\gamma_i}(x)\right)-\sum_{i=1}^n \beta_i \log\left(1-F_{1}^{\gamma_i}(x)\right),$ or equivalently, $\log \left(\frac{\prod_{i=1}^n \left(1-F_{2}^{\gamma_i}(x)\right)^{\beta_i}}{\prod_{i=1}^n \left(1-F_{1}^{\gamma_i}(x)\right)^{\beta_i}}\right)$ is increasing in $x$, which in turn gives    
\begin{equation}\label{e71}
\frac{\prod_{i=1}^n \left(1-F_{2}^{\gamma_{i}}(x)\right)^{\beta_i}}{\prod_{i=1}^n \left(1-F_{1}^{\gamma_{i}}(x)\right)^{\beta_i}}~is~increasing~in~x.
\end{equation}
Therefore,  from (\ref{e61}) and (\ref{e71}) it can be written that $$\left[\frac{\prod_{i=1}^n \left(1-F_{1}^{\gamma_{i}}(x)\right)^{\beta_i}}{\prod_{i=1}^n \left(1-F_{1}^{\alpha_{i}}(x)\right)^{\beta_i}}\right].\left[\frac{\prod_{i=1}^n \left(1-F_{2}^{\gamma_{i}}(x)\right)^{\beta_i}}{\prod_{i=1}^n \left(1-F_{1}^{\gamma_{i}}(x)\right)^{\beta_i}}\right]=\frac{\prod_{i=1}^n \left(1-F_{2}^{\gamma_{i}}(x)\right)^{\beta_i}}{\prod_{i=1}^n \left(1-F_{1}^{\alpha_{i}}(x)\right)^{\beta_i}}$$ is also increasing in $x$, implying that $U_{1:n}\le_{hr}V_{1:n}.$ \hfill$\Box$ \\
That the condition `for all $s>0$ $X_1^{s}\le_{hr}X_2^{s}$' of the previous theorem is only sufficient condition, is shown in the next counterexample. 
\begin{e1}\label{ce3}
Let $X_1$ and $X_2$ be two random variables having distribution functions $F_1(x)=1-e^{-3x^{4.4}}$ and $F_2(x)=1-e^{-0.2x^{0.4}}$ respectively. Now figure 4.1 (i) shows that although $X_1^{s}\le_{hr}X_2^{s}$ for $s=0.02$, figure 4.1 (ii) shows that there exists no hr ordering between $X_1^s$ and $X_2^s$ for $s=1.98$. Again, it can be shown that the same can be concluded for $s=0.01$ and $1.99$ respectively. Again, if $\mbox{\boldmath $\alpha$}=\left(1.99,0.01\right)\in \mathcal{D}_+,\mbox{\boldmath $\gamma$}=\left(1.98,0.02\right)\in \mathcal{D}_+$, $\mbox{\boldmath $\beta$}=\left(1,2\right)\in \mathcal{E}_+$ are taken, then Figure 4.2 shows that $ U_{1:n}\le_{hr}V_{1:n}.$ It is to be mentioned here that while plotting the curve the substitution $x=-\ln y$ has been used. 
\begin{figure}[ht]
\centering
\begin{minipage}[b]{0.45\linewidth}
\includegraphics[height=7 cm]{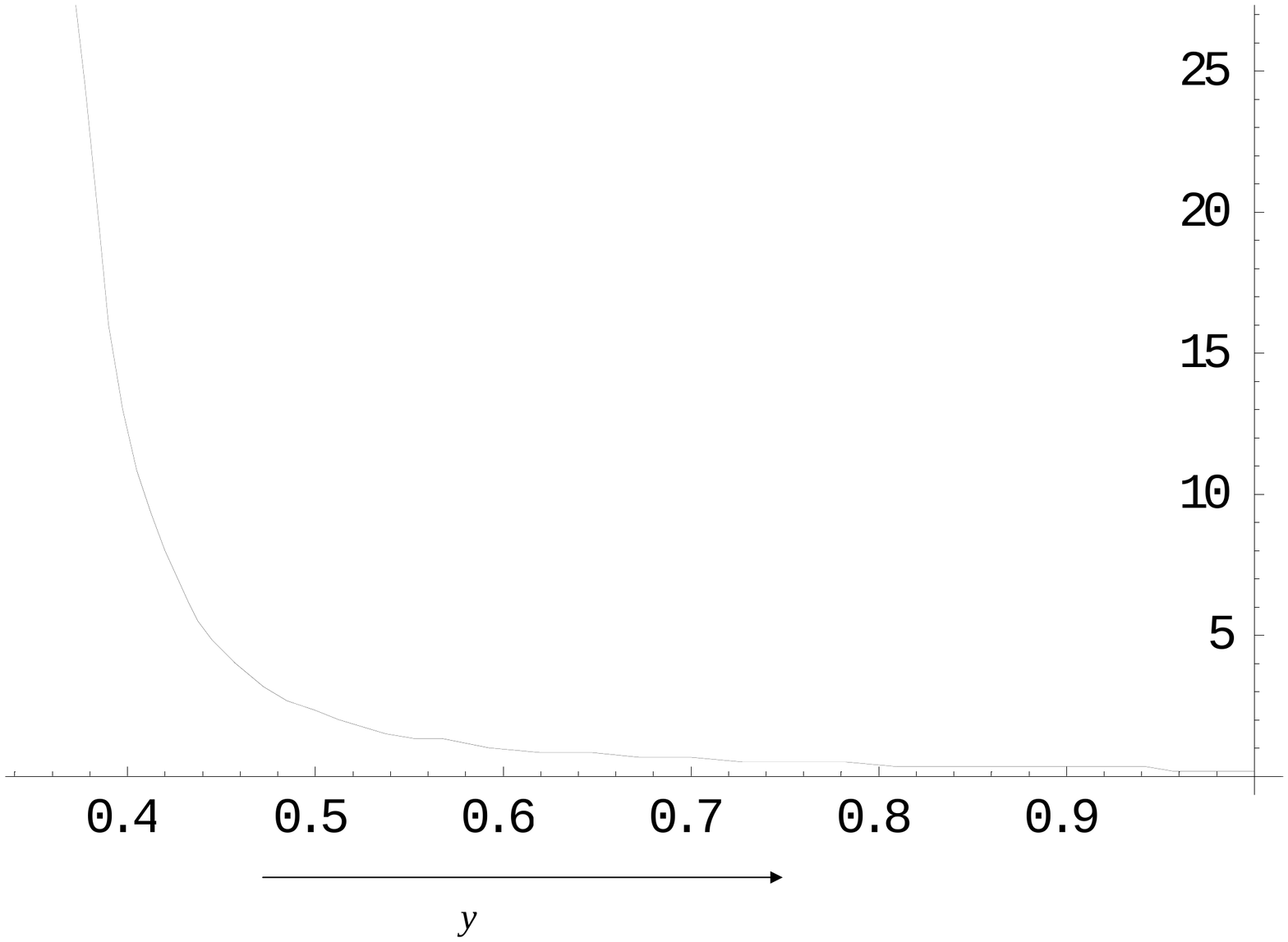}
\centering{$\left(i\right)$ Graph for $s=0.02$}
\end{minipage}
\quad
\begin{minipage}[b]{0.45\linewidth}
\includegraphics[height=7 cm]{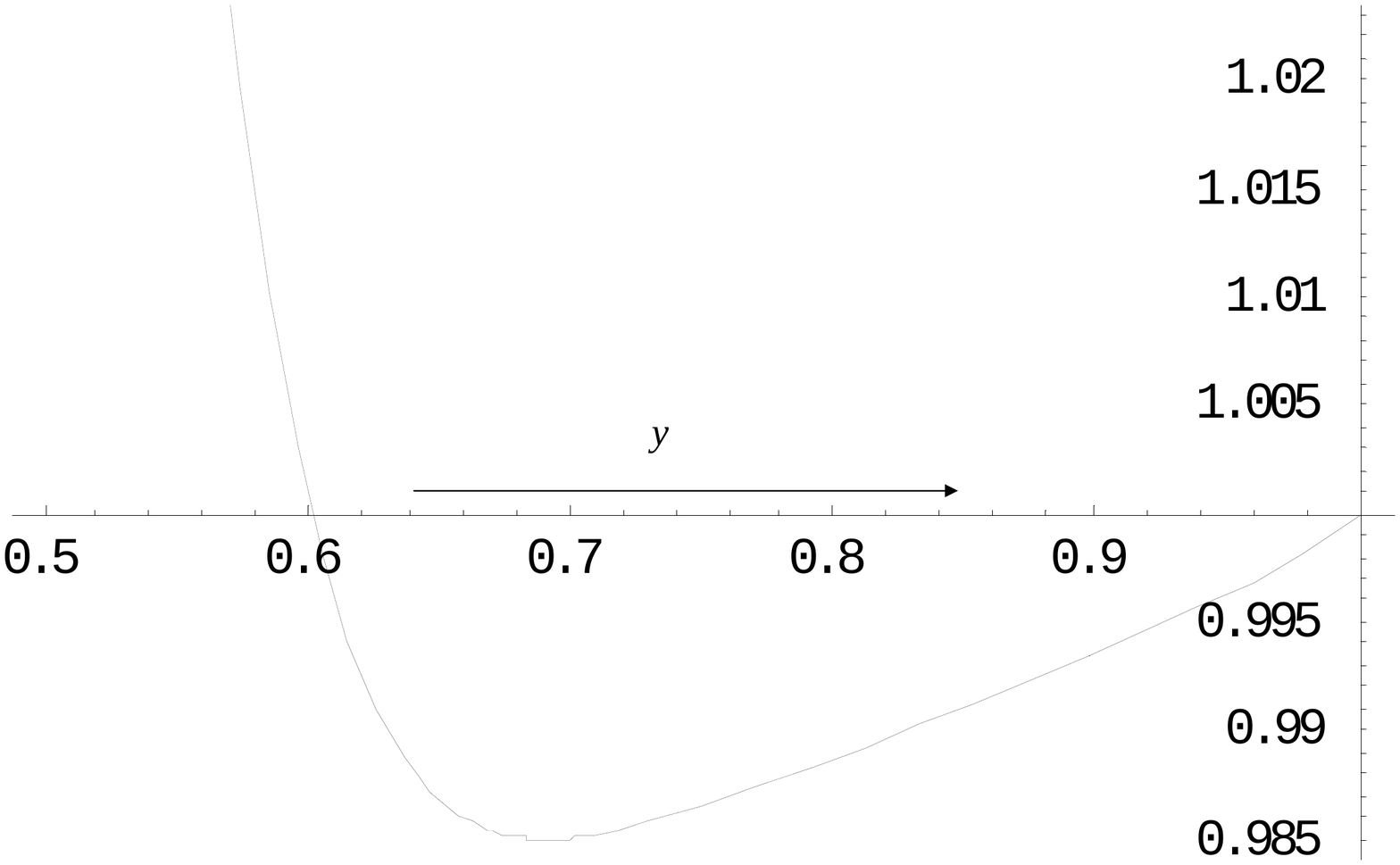}
\centering{$\left(ii\right)$ Graph for $s=1.98$}
\end{minipage}
\center{\it Figure 4.1: Graph of $\frac{1-F_{2}^s(x)}{1-F_{1}^s(x)}$ for $x=-\ln y$.} 
\end{figure}
\begin{figure}[t]\centering
\includegraphics[height=7 cm]{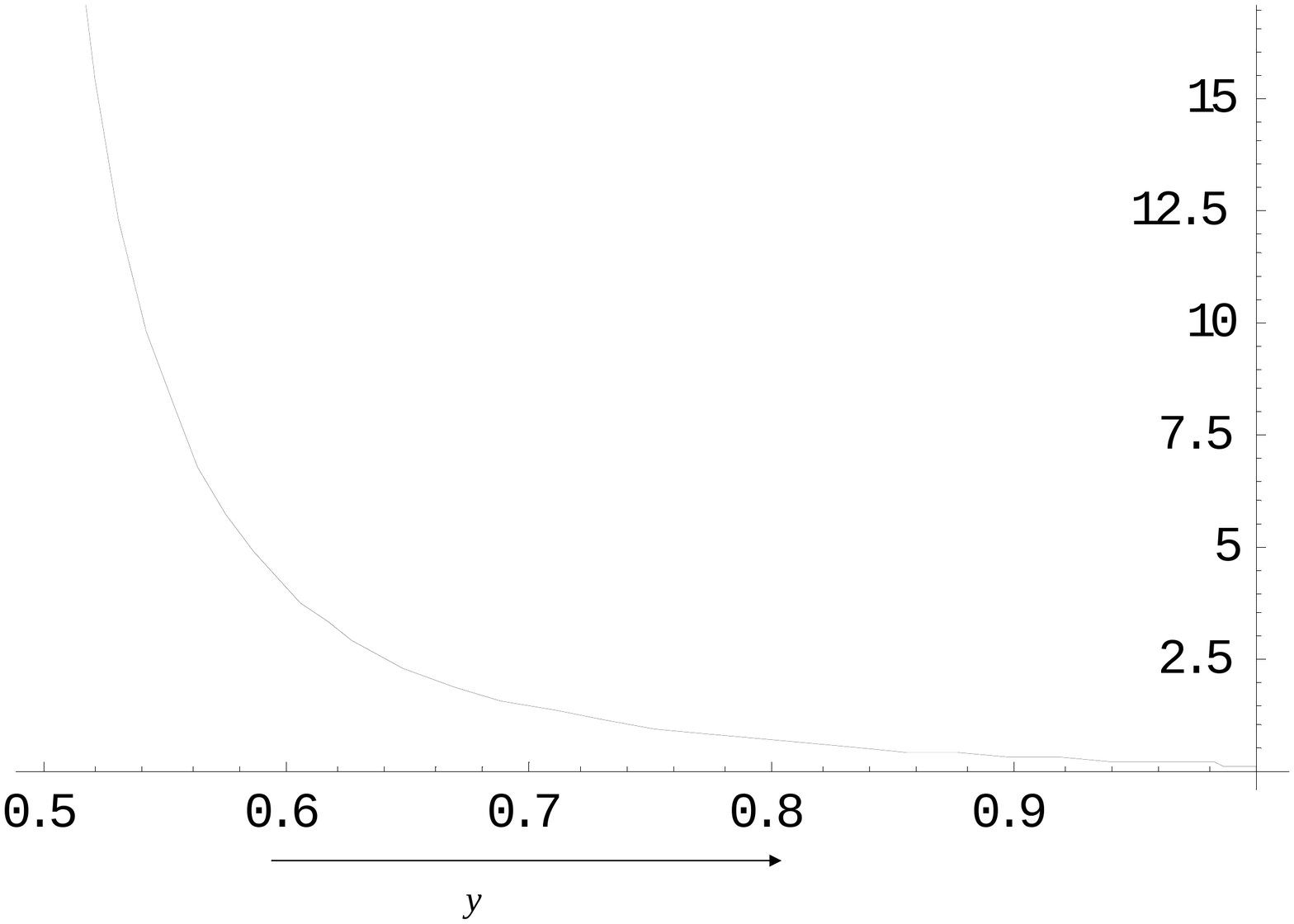}
\center{\it Figure 4.2: Graph of $\frac{\overline{H}_{1:n}\left(x\right)}{\overline{G}_{1:n}\left(x\right)}$ for $x=-\ln y$} 
\end{figure}
\end{e1}
\begin{t1}\label{th14}
Let $X_1$ and $X_2$ be two random variables having distribution functions $F_1$ and $F_2$ respectively. For $i=1,2,\ldots, n$, let $U_i$ and $V_i$ be two sets of mutually independent random variables with $U_i\sim $$Kw$-G$\left(\alpha_i,\beta_i,F_1\right)$ and $V_i\sim$ $Kw$-G$\left(\alpha_i,\delta_i,F_2\right)$. Further, for any real number $s>0$, suppose $X_i^{s}$ be a random variable having distribution function $F_i^{s}(x)$. If $\mbox{\boldmath $\beta$}\stackrel{m}{\succeq}\mbox{\boldmath $\delta$},$ $\mbox{\boldmath $\alpha$}\in \mathcal{D}_+$ and
\begin{enumerate}
\item[i)] $\mbox{\boldmath $\beta$},\mbox{\boldmath $\delta$}\in \mathcal{E}_+$, then $X_1\le_{hr}X_2$ implies $U_{1:n}\le_{hr}V_{1:n}$;
\item[ii)] $\mbox{\boldmath $\beta$},\mbox{\boldmath $\delta$}\in \mathcal{D}_+$, then $X_1\ge_{hr}X_2$ implies $U_{1:n}\ge_{hr}V_{1:n}.$
\end{enumerate}
\end{t1}
{\bf Proof:} Considering $W_i\sim $$Kw$-G$\left(\alpha_i,\delta_i,F_1\right)$, and using the same logic as of Theorem~\ref{th13} the theorem can be proved with the help of Theorem~\ref{th8}.\hfill$\Box$\\
\section{Concluding Remarks}
\hspace*{0.3 in} In this paper, we compare the hazard rate functions of the smallest order statistic arising from
independent heterogeneous $Kw$-G distributions when the shape parameters are majorized. The results are derived on the assumption that the parent cdf of the $Kw$-G random variables are homogenous and can be either identical or different. It is also shown that if the vectors of the shape parameters of the underlying distributions are in majorization order, then likelihood ratio ordering exists between the smallest order statistic from multiple-outlier $Kw$-G model with identical parent cdf $F.$ The results of this paper are applicable to a wide variety of distributions generated from $Kw$ distribution through the cdf $F$ as discussed in the Introduction, viz. $Kw$-N, $Kw$-W, $Kw$-Ga, $Kw$-Gu etc. 

\end{document}